\documentclass[preprint,authoryear]{elsarticle}
\usepackage{multirow}
\usepackage{color}
\usepackage{amsmath}
\usepackage{amssymb}
\usepackage{caption}
\usepackage{subcaption}
\usepackage{epstopdf}
\usepackage[colorinlistoftodos,textsize=small]{todonotes}
\usepackage[breaklinks]{hyperref}
\usepackage{tikz}
\usepackage{hypernat}

 \newcommand{\commentGeneric}[2]{\fbox{\begin{minipage}{0.6\textwidth}
 \textbf{#1:}~
 \begin{minipage}{0.8\textwidth}
 \it #2
 \end{minipage}\end{minipage}}\quad}
\newboolean{SHOWCOMMENTS}
\setboolean{SHOWCOMMENTS}{false}
\ifthenelse{\boolean{SHOWCOMMENTS}}{%
\newcommand{\commentU}[1]{\commentGeneric{Uli}{#1}}
}{
\newcommand{\commentU}[1]{}
}

\renewcommand{\vec}[1]{\mathbf{#1}}
\renewcommand{\a}[0]{\ensuremath{\alpha}}
\renewcommand{\b}[0]{\ensuremath{\beta}}
\renewcommand{\c}[0]{\ensuremath{\gamma}}

\newcommand{\e}[0]{\ensuremath{\epsilon}}
\newcommand{\x}[0]{\ensuremath{\vec{x}}}
\renewcommand{\u}[0]{\ensuremath{\vec{u}}}
\renewcommand{\c}[0]{\ensuremath{\vec{c}}}
\renewcommand{\l}[0]{\ensuremath{\lambda}}
\renewcommand{\L}[0]{\ensuremath{\Lambda}}

\graphicspath{ {./figures/} }

\begin{document}

\begin{frontmatter}
\title{Boundary Conditions for Free Interfaces with the Lattice Boltzmann Method}

\author[LSS]{Simon Bogner\corref{cor}}
\ead{simon.bogner@fau.de}

\author[LSS]{Regina Ammer}
\ead{regina.ammer@fau.de}

\author[LSS]{Ulrich R\"ude}
\ead{ruede@fau.de}

\address[LSS]{Lehrstuhl f\"ur Systemsimulation, 
	Universit\"at Erlangen-N\"urnberg,
	Cauerstra{\ss}e 11,
	91058 Erlangen,
        Germany}

 \cortext[cor]{Corresponding author}

\begin{abstract}
  In this paper we analyze the boundary treatment of the lattice Boltzmann method
(LBM) for simulating 3D flows with free surfaces. The widely used free surface
boundary condition of \citet{KoernerEtAl} is shown to be first order accurate.
The article presents a new free surface boundary scheme that is suitable for
second order accurate simulations based on the LBM. The new method takes into
account the free surface position and its orientation with respect to the
computational lattice. Numerical experiments confirm the theoretical findings
and illustrate the different behavior of the original method and the new
method.


\end{abstract}


\end{frontmatter}
\commentU{Hervorhebungen durch quotation marks vs. italic font - vereinheitlichen.}
\commentU{Unnötige Akronyme vermeiden ... habe wo es mir aufgefallen ist LHS RHS C:E. usw ersetzt}
\commentU{Gleichungen nur nummerieren, wo darauf Referenz genommen wird}

\section{\label{sec:introduction}Introduction}
\subsection{Motivation}
Since its early existence, the lattice Boltzmann method has been applied in
simulations of multi-phase flow phenomena
\citep{BenziEtAl,ChenDoolen98,AidunClausen}. For the major part, these efforts
have been based on diffusive interface theory
\citep{Do-QuangEtAl2000,Nourgaliev2003}. The \emph{free surface lattice
  Boltzmann method} (FSLBM) \citep{KoernerEtAl}, instead, is based on the
non-diffusive \emph{volume of fluid} (VOF) approach
\citep{Hirt68,ScardovelliZaleski99}, to track the motion of the interface and to
impose a free boundary condition locally. The dynamics of the gas phase is
neglected and a single-phase free boundary problem is solved instead of a
two-phase flow problem.
The method has been used successfully in the simulation of liquid-gas flows
of high viscosity and density ratio. Examples can be found in
\citet{Thuerey2004,Thuerey2006,XingEtAl2007a,XingEtAl2007b,Janssen2010,JanssenOffshore2010,Attar2011,DonathEtAl2010,Anderl2014,Ammer2014},
and in \citet{Bogner2013,Svec2012} for complex liquid-gas-solid flows.
However, no theoretical analysis of the FSLBM is currently available.  Hence,
the continued success in numerous applications motivates the interest in
developing a better theoretical foundation of the method.  In this paper we
present a detailed analysis of the free surface boundary condition as it is used
in the papers cited above. 
The analysis of lattice Boltzmann boundary schemes is mainly due to the works of
Ginzburg \citep{Ginzburg2003,Ginzburg2007,Ginzburg1994} and Junk
\citep{JunkYang2005,JunkYang2005b}. Here, we use a Chapman-Enskog ansatz similar
to \citet{Ginzburg2007}, to analyze the free surface boundary treatment.  We
show that the original FSLBM boundary condition, referred to as \emph{FSK} - rule
later on, is of first order in spatial accuracy. We then proceed to propose a
second order accurate free surface boundary condition as a possible
improvement. 
The new method is based on linear interpolation (\emph{FSL} - rule) and can be
analyzed by the same techniques mentioned above.

The considered numerical scheme including the free surface treatment is
introduced in Sec.~\ref{sec:method}. The analytic results, including the construction of the new
FSL - scheme, can be found in Sec.~\ref{sec:analysis}. We present various
numerical experiments in Sec.~\ref{sec:results} that confirm the predicted
behavior. Further discussion and outlook can be found in
Sec.~\ref{sec:conclusion}.

\subsection{Liquid Interfaces and Free Boundaries}
In this paragraph we introduce the model equations of a free surface.
Let here $\mu = \rho \nu$ denote the dynamic viscosity of the liquid.  The boundary
condition at a free interface with local unit normal $\vec{n}$ is given by
\begin{equation}
  (P - P_b + \sigma \kappa) n_\a = 2 \nu S_{\a \b} n_\b,
  \label{eq:freeboundary}
\end{equation}
where $P$ and $S_{\a \b} = \frac{1}{2} \left( \partial_{\a} j_{\b}
  + \partial_{\b} j_{\a} \right)$ are the local pressure and shear rate tensor
of the liquid, and $P_b$ is the surrounding pressure exerted on the
free boundary. The term $\sigma \kappa$ expresses the usual pressure jump due to
surface tension.  Let $\vec{\tau}$ be a unit vector tangential to the
interface. Projecting Eq.~(\ref{eq:freeboundary}) first on $\vec{n}$, and second, 
on $\vec{\tau}$ yields the conditions
\begin{subequations}
\label{eq:projectedfree}
 \begin{align}
  P - P_b + \sigma \kappa & = 2 \nu \left( \partial_n j_n \right), \label{eq:proj1} \\
  0 & = \partial_\tau j_n + \partial_n j_\tau \label{eq:proj2},
 \end{align}
\end{subequations}
for the normal and tangential viscous stresses, respectively. Here, $\partial_n$
and $\partial_{\tau}$ are the directional derivatives along $\vec{n}$ and
$\vec{\tau}$, respectively, while $j_n = n_{\a} j_{\a}$ and $j_{\tau} =
\tau_{\a} j_{\a}$. 

The remaining part of this paper deals with the construction
of lattice Boltzmann boundary rules, satisfying the above equations. The FSLBM
according to \citet{KoernerEtAl} uses a ``link-wise'' construction, as described in
\citet{Ginzburg2003,Ginzburg2007} for various Dirichlet- and mixed-type boundary 
conditions. 
We follow these works with our notation, such that the new boundary rules can 
be easily related to the ``multi-reflection'' context from the same authors. 
We remark, however, that there exist alternative, non - ``link-wise'' techniques to model free surface
boundary conditions \citep{GinzburgSteiner2003}, which are not considered in this article. 


\section{\label{sec:method}Numerical Method}
In this work, we 
will develop the FSLBM based on a \emph{two-relaxation-time} (TRT) collision
operator. This collision operator uses separate relaxation times for even and odd parts of the distribution function. 
This is particularly important when working with boundary conditions.  A
generalization to other collision operators is possible, but will not change
the convergence orders. The widespread \emph{single-relaxation-time} (SRT) collision
model, also known as lattice BGK model \citep{QianEtAl1992}, is a special case of
the TRT model, in which both relaxation times are chosen equal. Hence, the
results obtained in the analysis of Sec.~\ref{sec:analysis} can be transferred
directly to the SRT model.

\subsection{\label{sec:trt}Hydrodynamic TRT-model}
We assume a lattice Boltzmann equation \citep{BenziEtAl,ChenDoolen98,AidunClausen} with
two relaxation times according to \citet{Ginzburg2007,Ginzburg2005}. The evolution of the
distribution function $\vec{f} = (f_0, f_1, .., f_{Q-1})$ on the lattice for the finite
set of \emph{lattice velocities} $\{ \vec{c}_q \, | \, q=0, .., Q-1 \}$ is then described
by the equations
\begin{subequations}
\begin{align}
  f_q (\vec{x}+\vec{c}_q, t+1) &= \tilde{f_q} (\vec{x}, t),\label{eq:lbe1} \\
  \tilde{f}_q (\vec{x}, t) &= f_q (\vec{x}, t) + \lambda_{+}n_q^{+} + \lambda_{-}n_q^{-} + F_q,
  \label{eq:lbe2}
\end{align}
\end{subequations}
with $n_q^{\pm} = f_q^{\pm} - e_q^{\pm}$.
Here, Eq.~(\ref{eq:lbe1}) is referred to as the \emph{stream step}, and
Eq.~(\ref{eq:lbe2}) is the \emph{collision step}, yielding the \emph{post-collision}
distributions $\tilde{\vec{f}}$. The equation has two independent relaxation times
$\lambda_{+}, \lambda_{-} \in (-2, 0)$ for the \emph{even (symmetric)} and \emph{odd
  (anti-symmetric)} parts of the distribution function. $F_q$ is a source term that will
be discussed later. For the discrete range of values $q \in \{0, .., Q-1\}$, the opposite
index $\bar{q}$ is defined by the equation $\vec{c}_q = - \vec{c}_{\bar{q}}$ and thus,
\begin{align}
  f_q^{+} = \frac{1}{2}(f_q + f_{\bar{q}}),\;\text{ and }\; f_q^{-} = \frac{1}{2}(f_q - f_{\bar{q}}),
\end{align}
respectively. The equilibrium function $\vec{f}^{eq} = \vec{e}^{+} + \vec{e}^{-}$ is
of polynomial type with even part
\begin{equation}
  e_{q}^{+}(\rho, \vec{u}) = \frac{w_{q}}{c_s^2} \Pi_q,
  \label{eq:eqEven}
\end{equation}
where $\Pi_q = P + N_q$, with pressure $P$ and the non-linear contribution
\begin{equation}
  N_q = \frac{1}{2} \rho_0 u_{\a} u_{\b} \left( \frac{c_{q, \a} c_{q, \b}}{c_s^2} - \delta_{\a \b}\right),
\end{equation}
and odd equilibrium component
\begin{equation}
  e_{q}^{-}(\rho, \vec{u}^{eq}) = \frac{w_{q}}{c_s^2} \rho_0 c_{q, \a} u^{eq}_{\a}.
  \label{eq:eqOdd}
\end{equation}
Hereby, the pressure is defined 
by $P = c_s^2 \rho$.
The popular ``compressible''
form used in \cite{QianEtAl1992} is obtained if one sets $\rho_0=\rho$.
For $\rho_0 = 1$,
the ``incompressible'' equilibrium of \cite{HeLuo1997b} is obtained.
The \emph{lattice weights} $w_q = w_{|\vec{c}_q|}$ are chosen as in \citet{QianEtAl1992}, with $c_s= 1/\sqrt{3}$ as the corresponding \emph{lattice speed of sound}.
The non-linear part $N_q$ can be dropped for the simulation of Stokes-flow.
Macroscopic quantities are defined as moments of $\vec{f}$. In particular, the moments of
zeroth and first order,
\begin{align} 
  \rho = \frac{1}{c_s^2} P &= \sum_q f_q,
  \label{eq:density} \\
  \rho_0 \vec{U} = \rho_0 \vec{u} - \frac{\vec{F}}{2} &= \sum_q \vec{c}_q f_q,
  \label{eq:momentum}
\end{align}
define the pressure $P$ and fluid velocity $\vec{u}$. The shift in the fluid momentum by
$\vec{F}/{2}$ is necessary if external forces such as gravitation are included in
simulations \citep[cf.][]{BuickGreated2000,Ginzburg2007}.
For the latter, one can either
make use of additional force terms $F_q$ as in \citet{Guo2002,BuickGreated2000} and set
$\vec{u}^{eq}=\vec{U}$, or equivalently work without the source term $F_q$ and use
$\vec{u}^{eq}=\vec{U}-\vec{F}/{\lambda_{-}}$ instead. 
The latter simplifies the analysis
\citep{Ginzburg2007} and is adopted in the following.
The fluid momentum is $\vec{j} = \rho_0 \vec{u}$.
The lattice viscosity of the model is related to the symmetric relaxation time, via
\begin{equation}
  \nu = -\frac{1}{3}(\frac{1}{\l_+} + \frac{1}{2}).
\end{equation}

\subsection{Boundary Conditions}
The node positions $\x$ are restricted to a discrete subset (lattice) of nodes within a
bounded domain $\Omega \subset \mathbb{R}^n$. 
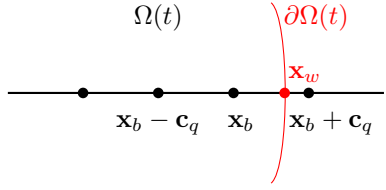
\begin{figure}
  \centering
  \begin{tikzpicture}
   \draw[thick] (0,0)--(5,0);
   \fill[black] (1,0) circle (2pt);
   \fill[black] (2,0) circle (2pt);
   \fill[black] (3,0) circle (2pt);
   \fill[black] (4,0) circle (2pt);
   \fill[red] (3.68,0) circle (2pt);
   \draw[red] (3.5,1.2) .. controls (3.75,1) and (3.75,-1) .. (3.5,-1.2);
   \node at (2,+1) {$\Omega(t)$};
   \node at (3.1,-0.4) {$\vec{x}_{b}$};
   \node at (2,-0.4) {$\vec{x}_{b}-\vec{c}_{q}$};
   \node at (4.3,-0.4) {$\vec{x}_{b}+\vec{c}_{q}$};
   \node[red] at (3.95,+0.25) {$\vec{x}_{w}$};
   \node[red] at (4.1,1) {$\partial \Omega(t)$};
  \end{tikzpicture}
  \caption{Schematic view of boundary node $\vec{x}_b \in \Omega(t)$ and boundary intersection point $\vec{x}_w \in \partial \Omega(t)$ along a lattice link $\vec{c}_q$.}
  \label{fig:boundary_node}
\end{figure}
A node $\x$ is called \emph{boundary node},
if the set of boundary links, $L_b(\x) := \{q \, | \, \x + \c_q \notin \Omega \}$, is
nonempty (cf. Fig.~\ref{fig:boundary_node}). If $\x_b$ is a boundary node, then for each $q \in L_b(\x_b)$, there is an
intersection $\x_w = \x_b + \delta \c_q \in \partial \Omega$ with $0 \leq \delta \leq 1$.
The value $f_{\bar{q}}(\x, t+1)$ then cannot be computed from
Eq.s~(\ref{eq:lbe1}-\ref{eq:lbe2}), and must be given in the form of a closure
relation. For this paper we consider linear link-wise closure relations 
that take the general form
\begin{equation}
  f_{\bar{q}}(\x_b, t+1) = a_0 \tilde{f}_{q}(\x_b, t) + \bar{a}_0 \tilde{f}_{\bar{q}}(\x_b, t) + a_1 \tilde{f}_{q}(\x_b - \c_q, t)\\ + f_q^{p.c.}(\x_b, t) + f^b_q(\x_w, t).
  \label{eq:closure}
\end{equation}
This can be categorized as a linear \emph{multi-reflection} closure rule
\citep{Ginzburg2003,Ginzburg2007}, with $\kappa_1=a_0$, $\bar{\kappa}_{-1}=\bar{a}_0$,
$\kappa_0=a_1$, and $\kappa_{-1}=\kappa_{-2}=0$ when using the notation of the respective
articles.
$f_q^{p.c}$ is a term depending on the local non-equilibrium component
$n_q(\x_b,t)$, and $f^b_q(\x_w, t)$ 
depends on the (macroscopic) boundary values at the wall point $\x_w$.
For Dirichlet-type boundary conditions on the pressure or the momentum, respectively, 
this term takes the form
\begin{equation}
  f^b_q = \alpha_+ e_q^+(\rho_b, \u_b) + \alpha_- e_q^-(\rho_b, \u_b),
\end{equation}
where the $\alpha_{\pm}$ are linear combinations of the coefficients $a_0$, $\bar{a}_0$,
$a_1$, depending on the specific boundary condition.

\subsection{\label{sec:fslbm}Free Surface Lattice Boltzmann Method}
For the FSLBM according to \citet{KoernerEtAl}, the lattice Boltzmann
equation scheme described above is extended by a \emph{volume-of-fluid} indicator
function $\varphi(\vec{x}, t) \in [0, 1]$ \citep{HirtNichols,ScardovelliZaleski99}. This
function is defined as the volume fraction of liquid within the cubic unit cell centered
around the lattice node at $\vec{x}$, thus giving an implicit description of the free
surface between liquid and gas.
For dynamic simulations the indicator function $\varphi$
must be advected after each time step.
It represents a boundary for the hydrodynamic
simulation, and its closure relation as given in \citet{KoernerEtAl} reads
\begin{equation}
  f_{\bar{q}}(\x_b, t+1) = - \tilde{f}_q(\x_b,t) + 2 \cdot e_q^{+}(\rho_b, \u_b),
  \label{eq:koerner}
\end{equation}
where $P_b = c_s^2 \rho_b$ is the boundary value for the pressure at the free
surface, and $\u_b$ is the velocity of the interface and must be extrapolated to
the boundary from the nodes.  The lattice Boltzmann domain $\Omega(t)$ is thus
limited to nodes $\x$ with $\varphi(\x,t)>0$. 
Eq.~(\ref{eq:koerner}) is applied at \emph{interface nodes} $\x_b$ for all links
$q$ that are connected to inactive \emph{gas nodes} with $\varphi(\x_b+\c_q, t) =
0$. Active lattice Boltzmann nodes in $\Omega(t)$ are also called \emph{liquid
 nodes}.

Surface tension can be directly incorporated in the FSLBM by including also the
Laplace pressure term $\sigma \kappa$ in Eq.~\eqref{eq:koerner} in place of
$P_b$. As other interface capturing methods, this requires a local approximation
of the interface curvature $\kappa$ \citep{ScardovelliZaleski99,Fuster2009}.




\section{\label{sec:analysis}Free surface Boundary Conditions}
\subsection{Chapman-Enskog Analysis}
We apply the Chapman-Enskog ansatz of \citet{Ginzburg2007} for
incompressible flow.
Based on diffusive time scaling \citep{JunkYang2005}, the time-derivatives of the
first order in the expansion parameter $\e$ are dropped and one seeks solutions to the
system of Eq.s~(\ref{eq:lbe1}-\ref{eq:lbe2}) with the scaled space and time step satisfying
$\Delta x^2 = \Delta t = \mathcal{O}(\e^2)$.
For brevity, we introduce $\partial_q :=
c_{q,\a} \partial_{\a}$, $j_q := c_{q,\a} j_{\a}$.
Then, the non-equilibrium solution up
to the third order, split into even and odd parts, for constant external forcing reads
\begin{equation}
  n_q^{\pm} = \frac{1}{\l_{\pm}} \left[ \partial_q ( e_q^{\mp} - \L_{\mp} \partial_q e_q^{\pm}  ) + \partial_t e_q^{\pm} \right] + \mathcal{O}(\e^3),
  \label{eq:ce}
\end{equation}
where the coefficients $\L_{\mp}$ are defined as
\begin{equation}
  \L_{\pm} = -(\frac{1}{2} + \frac{1}{\l_{\pm}}).
\end{equation}
We further define the product $\L := \L_+ \L_-$, which is useful to characterize the
parametrization of the model. Substituting the polynomial equilibrium,
Eq.s~(\ref{eq:eqEven}-\ref{eq:eqOdd}) and considering only constant external forcing, we
can directly express $n_q$ in terms of macroscopic variables by
\begin{subequations}
\begin{align}
  n_q^{+} &= \frac{1}{\l_+} \frac{w_q}{c_s^2} \left[ \partial_q \left(j_q - \L_- \partial_q \Pi_q \right) + \partial_t \Pi_q \right], \label{eq:nqp} \\
  n_q^{-} &= \frac{1}{\l_-} \frac{w_q}{c_s^2} \left[ \partial_q (\Pi_q -
    \L_+ \partial_q j_q) +\partial_t j_q \right]. \label{eq:nqm}
\end{align}
\end{subequations}
The approximate solution based on Eq.~(\ref{eq:ce}) can be used to analyze boundary
conditions, by substituting into the respective closure relation and rewriting it for the
macroscopic variables in question, after Taylor-expanding all occurrences of $f_q$ around
$(\x_b,t)$. Notice that space and time derivatives are of first and second order in
$\epsilon$, respectively, and only terms up to $\mathcal{O}(\epsilon^2)$ need to be
included in the analysis. Based on Eq.~(\ref{eq:ce}), it is possible to construct new
boundary schemes by substituting into a general form like Eq.~(\ref{eq:closure}) and then
matching the unknown coefficients to yield the desired condition for the macroscopic
variables. We will apply this technique to derive a higher order free surface boundary
condition in Sec.~\ref{sec:fsl}. We recall that the analysis of the present paper holds
for the popular SRT collision model with only one relaxation time $\l_+ = \l_-$, too, if
all occurrences of the anti-symmetric relaxation time $\l_-$ are replaced with the
symmetric relaxation parameter $\l_+$, that controls the viscosity
(cf. Sec.~\ref{sec:method}).

\subsection{Analysis of the FSLBM}
\label{sec:}
The free surface boundary condition of Eq.~(\ref{eq:koerner}) is expanded around
$(\x_b,t)$ on the left hand side. Substituting $f_q = e_q + n_q$ to split equilibrium and
non-equilibrium parts of the solution, and further separating into even and odd parts, one
obtains up to the order $\e^2$,
\begin{equation}
  \left[e_q^{+} -\L_+ \l_+ n_q^+  + \frac{\l_-}{2} n_q^- + \partial_t(e_q^+ - e_q^-) \right](\x_b, t) = e_q^{+}(\x_w).
  \label{eq:expFSLBM}
\end{equation}
Notice, that all terms associated with the point $(\x_b, t)$ have been collected on the
left hand side, while on the right hand side only a boundary value term associated with
intersection point $\x_w$ remains. Substituting the second order non-equilibrium solution,
the left hand side of Eq.~(\ref{eq:expFSLBM}) results in
\begin{equation}
  \left[ \left( 1 + \frac{1}{2}\partial_q + \L\partial_q^2 +(1-\L_+)\partial_t \right) e_q^+ - \left( \L_+\partial_q + \frac{\L_+}{2}\partial_q^2 + \frac{1}{2}\partial_t \right) e_q^-  \right](\x_b,t).
\end{equation}
Finally, using the polynomial equilibria of
Eq.s~(\ref{eq:eqEven}-\ref{eq:eqOdd}), neglecting the non-linear terms and
dropping all time derivatives, we obtain
\begin{equation}
  \left[ (1 + \frac{1}{2}\partial_q + \L\partial_q^2 ) P - \L_+ (1 +\frac{1}{2}\partial_q ) \partial_q j_q  \right](\x_b, t) = P_b(\x_w,t).
  \label{eq:fsfit}
\end{equation}
Obviously, the left hand side of Eq.~(\ref{eq:fsfit}) can be interpreted as a combination of the
Taylor-series approximation of the pressure $P$ and shear rate $\partial_q j_q$
at the point $\x_b + 1/2 \c_q$.
Hence, assuming $\delta=1/2$, Eq.~(\ref{eq:fsfit}) implies a second order (third order, for $\L=1/8$) accurate
agreement of the pressure with boundary value $P_b$, combined with a second
order condition of vanishing shear stress in $\x_w$.
One can show analytically or by numerical experiment that in the special case of a
steady parabolic force-driven
tangential free-surface flow over a lattice aligned plane
with no-slip boundary condition
is solved without error by the FSLBM when the boundary condition of
Eq.~(\ref{eq:koerner}) is applied and if the film-thickness is an 
integer value such that $\delta = 1/2$ (cf. also Sec.~\ref{sec:filmflow}).
However, if $\delta \neq 1/2$, as in most relevant cases, then the
 spatial accuracy for both pressure and
shear drops to the first order.
Also, this boundary rule fulfills
Eq.~\eqref{eq:proj2} only, but does not include the normal viscous stress term of
Eq.~\eqref{eq:proj1}. We will show in Sec.~\ref{sec:fsl} how this can be improved.

\subsection{Second order boundary condition for the shear rate}
\label{sec:shearbc}
Starting from Eq.~(\ref{eq:closure}) it is possible to construct a higher order boundary
condition for pressure and shear stress. To this end, we use the local correction term
\begin{equation}
  f_q^{p.c.}(\x_b,t) = C \cdot n_q^+(\x_b,t),
  \label{eq:correction}
\end{equation}
and a boundary value term of the form
\begin{equation}
  f_q^{b} = \alpha_+ e_q^+(\rho_b, \u_b) + D \cdot c_{q,\a} c_{q,\b} S_{\a \b}^b,
  \label{eq:bterm}
\end{equation}
which allows to prescribe the boundary values, $P_b = c_s^2 \rho_b$ for pressure, and $S_{\a \b}^b$ for the
shear rates in $\x_w$.  
We use $\tilde{f}_q(\x_b-\c_q, t) = f_q(\x_b, t+1) \approx
f_q(\x_b, t) + \partial_t e_q(\x_b,t)$, and then rewrite Eq.~(\ref{eq:closure}) placing
all terms except the boundary value 
$f^b_q(\x_w,t)$ on the left hand side.
Using the Chapman-Enskog
approximation from Eq.~(\ref{eq:ce}) and rearranging terms, we obtain
\begin{equation}
  \left[ \alpha_+ e_q^+ + \beta_+ n_q^+ + \alpha_- e_q^- + \beta_- n_q^- + \alpha_+^t\partial_t e_q^+ +\alpha_-^t\partial_t e_q^- \right](\x_b,t) = f_q^b(\x_w),
  \label{eq:fitting}
\end{equation}
where
\begin{subequations}
\begin{align}
  \alpha_+ &= 1 - a_0 -\bar{a}_0 - a_1, \label{eq:cl1} \\
  \beta_+ &= 1 - (1+\l_+)(a_0 + \bar{a}_0) - a_1 - C, \\
  \alpha_- &= \bar{a}_0 - a_0 - a_1 -1, \\ 
  \beta_- &= (1+\l_-) \bar{a}_0 -(1+\l_-) a_0 -a_1 -1, \\
  \alpha_+^t &= 1-a_1, \\
  \alpha_-^t &= -(1+a_1).
\end{align}
\end{subequations}
The aim of the following construction is to match these coefficients with the spatial
Taylor-series around $\x_b$ up to the second order for pressure and shear rate,
respectively. As the spatial derivatives of pressure and momentum are contained in the
non-equilibrium functions, Eq.s~(\ref{eq:nqp}-\ref{eq:nqm}), the system of equations
follows as
\begin{subequations}
\begin{align}
  \alpha_- &= 0, \label{eq:cl2} \\ 
  \beta_- &= \alpha_+ \delta \l_-, \label{eq:cl3} \\
  \beta_+ &= -\alpha_+ \L_+ \l_+, \label{eq:cl4}
\end{align}
\end{subequations}
keeping $\alpha_+$ as free parameter. Here, $\beta_-$ is chosen to fit the coefficient of
the first order derivative of the pressure in $n_q^-$. $\beta_+$ is chosen to fit the
coefficient of $\partial_q j_q$ from $n_q^+$ with the second order derivative
$\partial_q^2 j_q$ from $n_q^-$. The closure relation coefficients follow from the
Eq.s~(\ref{eq:cl1}, \ref{eq:cl2} and \ref{eq:cl3}) as
\begin{subequations}
  \label{eq:fsli}
\begin{align}
  a_0 &= 1 - \alpha_+(\frac{1}{2} + \delta), \label{eq:fsli1} \\
  \bar{a}_0 &= 1 - \frac{1}{2} \alpha_+, \label{eq:fsli2} \\
  a_1 &= \delta \alpha_+ -1, \label{eq:fsli3}
\end{align}
while the coefficient $C$ in the correction term $f_q^{p.c}$, as derived from Eq.~(\ref{eq:cl4}), is
\begin{equation}
  C = \alpha_+ \l_+ (\frac{1}{2} +\delta) - 2 \l_+. \label{eq:fsli4}
\end{equation}
\end{subequations}
Using now Eq.~(\ref{eq:ce}) to express Eq.~(\ref{eq:fitting}) in terms of gradients of the equilibria one obtains
\begin{equation}
  \begin{split}
  \alpha_+ \left( 1 + \delta \partial_q + \L \partial_q^2  \right)e_q^+ 
  - \alpha_+ \L_+ \left( 1 + \delta \partial_q \right) \partial_q e_q^- \\
  + \left( 1-a_1-\alpha_+ \L_+ \right) \partial_t e_q^+
  - \left( 1+a_1 -\alpha_+ \delta \right) \partial_t e_q^-
  = f_q^b,
  \end{split}
\end{equation}
and after substituting $a_1$ and $f^b_q$,
\begin{equation}
  \begin{split}
  \alpha_+ \left( 1 + \delta \partial_q + \L \partial_q^2  \right)e_q^+ 
  - \alpha_+ \L_+ \left( 1 + \delta \partial_q \right) \partial_q e_q^-
  + \left( 2 -\alpha_+ (\delta +\L_+) \right) \partial_t e_q^+ \\
  = \alpha_+ e_q^+(\rho_b, \u_b) -\alpha_+ \L_+ \frac{w_q}{c_s^2} c_{q,\a} c_{q,\b}
  S_{\a \b}^b.
  \label{eq:sli}
  \end{split}
\end{equation}
On the right hand side, the unknown coefficient of the
boundary term, Eq.~(\ref{eq:bterm}) is determined as $D = -\alpha_+ \L_+ \frac{w_q}{c_s^2}$ 
to fit with the left hand side.
Because the spatial approximation fits up to the second order for both pressure and shear rate, the
boundary condition can be classified of second order in space for both pressure and
momentum.  

Since $e_q^+$ contains the non-linear terms that are often responsible for
numerical instabilities, the corresponding error terms deserve special
attention: the spatial error of second order $\alpha_+ (\delta^2/2 -
\L)\partial_q e_q^+$ is bounded and independent of the viscosity if $\L$ is
fixed to a constant value, which is a usual requirement for parametrizations of
the TRT collision model \citep{Ginzburg2007}.  If the SRT model is used, then
$\L = (1/2 + 1/\l_+)^2$, because both relaxation times $\l_+$ and $\l_-$ are
identical. Hence, one should avoid values close to zero for $\l_+$, and exclude
very high lattice viscosities with the SRT model. The error in time, $\alpha_+(1
- {2}/{\alpha_+}+\delta+\L_+)\partial_t e_q^+$ depends through $\L_+/3 = \nu$ on
the lattice viscosity.  However, usually one either has high Mach numbers and
low viscosity (high Reynolds number regime) and hence $\L_+ \ll 1$, or a high
Mach number with high viscosities (low Reynolds number regime). The second case
arises typically if the LBM is used to simulate Stokes-like flow, and then the
non-linear terms in $e_q^+$ do not need to be included in the equilibrium
function \citep{Ladd1994a}. Hence, in this case the momentum-dependent error in
time may be eliminated.

It should be noted that the coefficients $a_0$, $\bar{a}_0$ and $a_1$ are the identical to
the linear interpolation based pressure boundary condition ``PLI'' of
\citet{Ginzburg2007}. This is a direct consequence of the construction described above of
matching the coefficients of the pressure gradients in the closure relation. The
coefficients $C$ and $D$, however, are different from the PLI - rule.
They are needed to obtain the $\partial_q j_q$ term in the left hand side
of Eq.~(\ref{eq:sli}), and to define the boundary value for the shear rate in the right 
hand side, respectively.

\subsection{Second order boundary condition for free surfaces}
\label{sec:fsl}
The boundary condition of the preceding Sec.~\ref{sec:shearbc} can be used to replace the
first order free surface rule of Eq.~(\ref{eq:koerner}).
In fact, the second order version
of Eq.~(\ref{eq:koerner}) is obtained by the defining Eq.s~(\ref{eq:fsli1}-\ref{eq:fsli4})
and setting $S_{\a \b}^b=0$ as boundary value.
However, for full consistency with the
physical model Eq.s~(\ref{eq:proj1}-\ref{eq:proj2}), it is necessary to control the
tangential and normal shear stresses individually.
Let $\{ \vec{t}_1, \vec{t}_2, \vec{n} \}$ be a local orthonormal
basis with $\vec{n}$ normal to the free boundary.
Using the
indices $\{ \a', \b' \}$ for the corresponding coordinate system, related to the standard
coordinates by rotation $l_{\a \b}$, the shear rate tensor can be expressed in the local
basis via 
\begin{equation}
  \bar{S}_{\a' \b'} = l_{\a \a'} l_{\b \b'} S_{\a \b}.
  \label{eq:bstress}
\end{equation}
The respective entries of the shear tensor $\bar{S}_{\a' \b'}$ can now be set
individually according to Eq.s~(\ref{eq:proj1}-\ref{eq:proj2}), leaving the
remaining components untouched.
In practice, $S_{\a \b}$ must be obtained by extrapolation from the bulk.

\begin{table}
  \centering
  \begin{tabular}{l|*{3}{c}cc|c}
          & $a_0$ & $\bar{a}_0$ & $a_1$ & $\alpha_+$ & $C$   & $D$ \\
    \hline
    FSK  & $-1$ & $0$ & $0$ & $2$ & $0$ & $-2 \Lambda_+ \frac{w_q}{c_s^2}$   \\
    FSL               & $\frac{1}{2}-\delta$ & $\frac{1}{2}$ & $\delta -1$ & 1 &  $\lambda_+ (\frac{1}{2} +\delta) - 2\lambda_+ $ & $-\Lambda_+ \frac{w_q}{c_s^2}$   \\
  \end{tabular}
  \caption{Coefficients of closure relation, Eq.~\eqref{eq:closure} with correction term $f_q^{p.c.}$ from Eq.~\eqref{eq:correction} and boundary value term $f_q^b$ from Eq.~\eqref{eq:bterm}.
    The FSK boundary rule is first order and purely local. FSL is second order and based on linear interpolation of the PDFs.\\
  }
  \label{tab:fsbc}
\end{table}

In Tab.~\ref{tab:fsbc} we have collected the coefficients for all the free surface
conditions considered in this paper. The FSK - rule is only first order except for a
plane aligned interface at distance $\delta_x/2$ from the boundary nodes, and equivalent
to the original FSLBM closure relation, Eq.~\eqref{eq:koerner}, if $D=0$. The
\emph{FSL}-rule is the free surface condition based on the construction of
Sec.~\ref{sec:shearbc}. This free surface condition is of second order spatial accuracy,
and fully consistent with Eq.s~(\ref{eq:proj1}-\ref{eq:proj2}). It should be noted, that
by setting $D=0$, we obtain simplified boundary conditions, consistent with Eq.~\eqref{eq:proj2} only, but
neglecting the normal viscous stresses in Eq.~\eqref{eq:proj1}. The importance of these
terms has been discussed for instance in \citet{McKibben95,Hirt68} and depends on the
respective problem. In fact, for $D=0$ all shear stress components vanish at the
boundary. Numerical simulations of free surface flows often use this simplified free
surface condition. In this case the $S_{\a \b}^b$ in Eq.~\eqref{eq:bterm} drops out, and
the condition can be implemented without the construction above and without extrapolation of $S_{\a \b}$.


\section{\label{sec:results}Numerical Results}
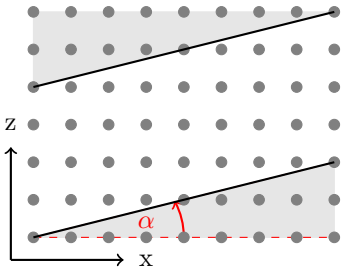
\begin{figure}
\centering
 \begin{tikzpicture}
  \filldraw[black!10](1,0)--(5,0)--(5,1)--(1,0);
  \filldraw[black!10](1,2)--(5,3)--(1,3)--(1,2);
  \draw[thick,->] (0.7,-0.3)--(2.2,-0.3);
  \draw[thick,->] (0.7,-0.3)--(0.7,1.2);
  \node at (2.5,-0.3) {x};
  \node at (0.7,1.5) {z};
  \draw[dashed,red] (1,0)--(5,0);
  \draw[dashed,red] (1,0)--(5,1);
  \draw[red,thick,->] (3,0) arc (0:29:1);
  \node[red] at (2.5,0.2) {$\alpha$};
  \filldraw[black!50] (1.0,0) circle (2pt);
  \filldraw[black!50] (1.5,0) circle (2pt);
  \filldraw[black!50] (2.0,0) circle (2pt);
  \filldraw[black!50] (2.5,0) circle (2pt);
  \filldraw[black!50] (3.0,0) circle (2pt);
  \filldraw[black!50] (3.5,0) circle (2pt);
  \filldraw[black!50] (4.0,0) circle (2pt);
  \filldraw[black!50] (4.5,0) circle (2pt);
  \filldraw[black!50] (5.0,0) circle (2pt);
  \filldraw[black!50] (1.0,0.5) circle (2pt);
  \filldraw[black!50] (1.5,0.5) circle (2pt);
  \filldraw[black!50] (2.0,0.5) circle (2pt);
  \filldraw[black!50] (2.5,0.5) circle (2pt);
  \filldraw[black!50] (3.0,0.5) circle (2pt);
  \filldraw[black!50] (3.5,0.5) circle (2pt);
  \filldraw[black!50] (4.0,0.5) circle (2pt);
  \filldraw[black!50] (4.5,0.5) circle (2pt);
  \filldraw[black!50] (5.0,0.5) circle (2pt);
  \filldraw[black!50] (1.0,1.0) circle (2pt);
  \filldraw[black!50] (1.5,1.0) circle (2pt);
  \filldraw[black!50] (2.0,1.0) circle (2pt);
  \filldraw[black!50] (2.5,1.0) circle (2pt);
  \filldraw[black!50] (3.0,1.0) circle (2pt);
  \filldraw[black!50] (3.5,1.0) circle (2pt);
  \filldraw[black!50] (4.0,1.0) circle (2pt);
  \filldraw[black!50] (4.5,1.0) circle (2pt);
  \filldraw[black!50] (5.0,1.0) circle (2pt);
  \filldraw[black!50] (1.0,1.5) circle (2pt);
  \filldraw[black!50] (1.5,1.5) circle (2pt);
  \filldraw[black!50] (2.0,1.5) circle (2pt);
  \filldraw[black!50] (2.5,1.5) circle (2pt);
  \filldraw[black!50] (3.0,1.5) circle (2pt);
  \filldraw[black!50] (3.5,1.5) circle (2pt);
  \filldraw[black!50] (4.0,1.5) circle (2pt);
  \filldraw[black!50] (4.5,1.5) circle (2pt);
  \filldraw[black!50] (5.0,1.5) circle (2pt);
  \filldraw[black!50] (1.0,2.0) circle (2pt);
  \filldraw[black!50] (1.5,2.0) circle (2pt);
  \filldraw[black!50] (2.0,2.0) circle (2pt);
  \filldraw[black!50] (2.5,2.0) circle (2pt);
  \filldraw[black!50] (3.0,2.0) circle (2pt);
  \filldraw[black!50] (3.5,2.0) circle (2pt);
  \filldraw[black!50] (4.0,2.0) circle (2pt);
  \filldraw[black!50] (4.5,2.0) circle (2pt);
  \filldraw[black!50] (5.0,2.0) circle (2pt);
  \filldraw[black!50] (1.0,2.5) circle (2pt);
  \filldraw[black!50] (1.5,2.5) circle (2pt);
  \filldraw[black!50] (2.0,2.5) circle (2pt);
  \filldraw[black!50] (2.5,2.5) circle (2pt);
  \filldraw[black!50] (3.0,2.5) circle (2pt);
  \filldraw[black!50] (3.5,2.5) circle (2pt);
  \filldraw[black!50] (4.0,2.5) circle (2pt);
  \filldraw[black!50] (4.5,2.5) circle (2pt);
  \filldraw[black!50] (5.0,2.5) circle (2pt);
  \filldraw[black!50] (1.0,3.0) circle (2pt);
  \filldraw[black!50] (1.5,3.0) circle (2pt);
  \filldraw[black!50] (2.0,3.0) circle (2pt);
  \filldraw[black!50] (2.5,3.0) circle (2pt);
  \filldraw[black!50] (3.0,3.0) circle (2pt);
  \filldraw[black!50] (3.5,3.0) circle (2pt);
  \filldraw[black!50] (4.0,3.0) circle (2pt);
  \filldraw[black!50] (4.5,3.0) circle (2pt);
  \filldraw[black!50] (5.0,3.0) circle (2pt);
  \draw[thick] (1,0)--(5,1);
  \draw[thick] (1,2)--(5,3);
 \end{tikzpicture}
 \caption{Sketch of channel flow, rotated by angle $\alpha$ with respect to the
   lattice.}
  \label{fig:sketch_channelFlow}
\end{figure}

All test cases presented in the following have been conducted using the TRT
collision operator described in Sec.~\ref{sec:trt} with a $D3Q19$ lattice model.
The sketch of numerical test cases for the different channel flows is depicted
in Fig.~\ref{fig:sketch_channelFlow}. Hereby, the flow variables are assumed
constant along the $y$-axis, and the channel is rotated by an angle $\alpha$
about the $y$-axis. The test cases have been realized using the waLBerla
\citep{walberla2011,Thürey2006} framework.

\subsection{Transient Evolution of Plate-Driven Planar Flow}
\begin{figure}
    \centering
    \includegraphics[width=0.66\textwidth]{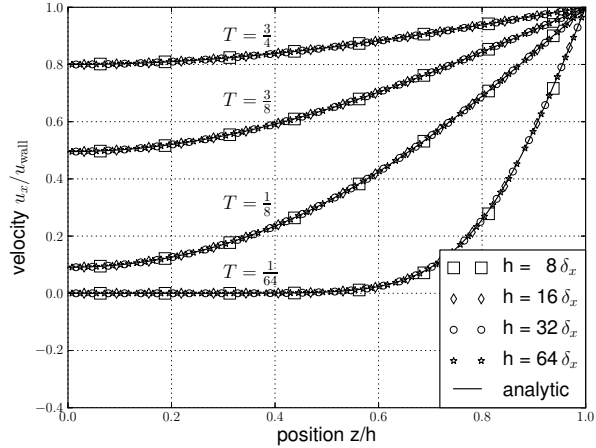}
    \caption{The velocity profile at non-dimensional times $T=1/64$, $1/8$,
      $3/8$. Both the original free surface boundary condition (FSK) and the new
      boundary condition based on linear interpolation (FSL) are very close to
      the analytical formula, Eq.~\eqref{eq:freeUniform}.}
    \label{fig:freeUniformProfile}
\end{figure}

\begin{figure}
  \centering
  \includegraphics[width=0.49\linewidth]{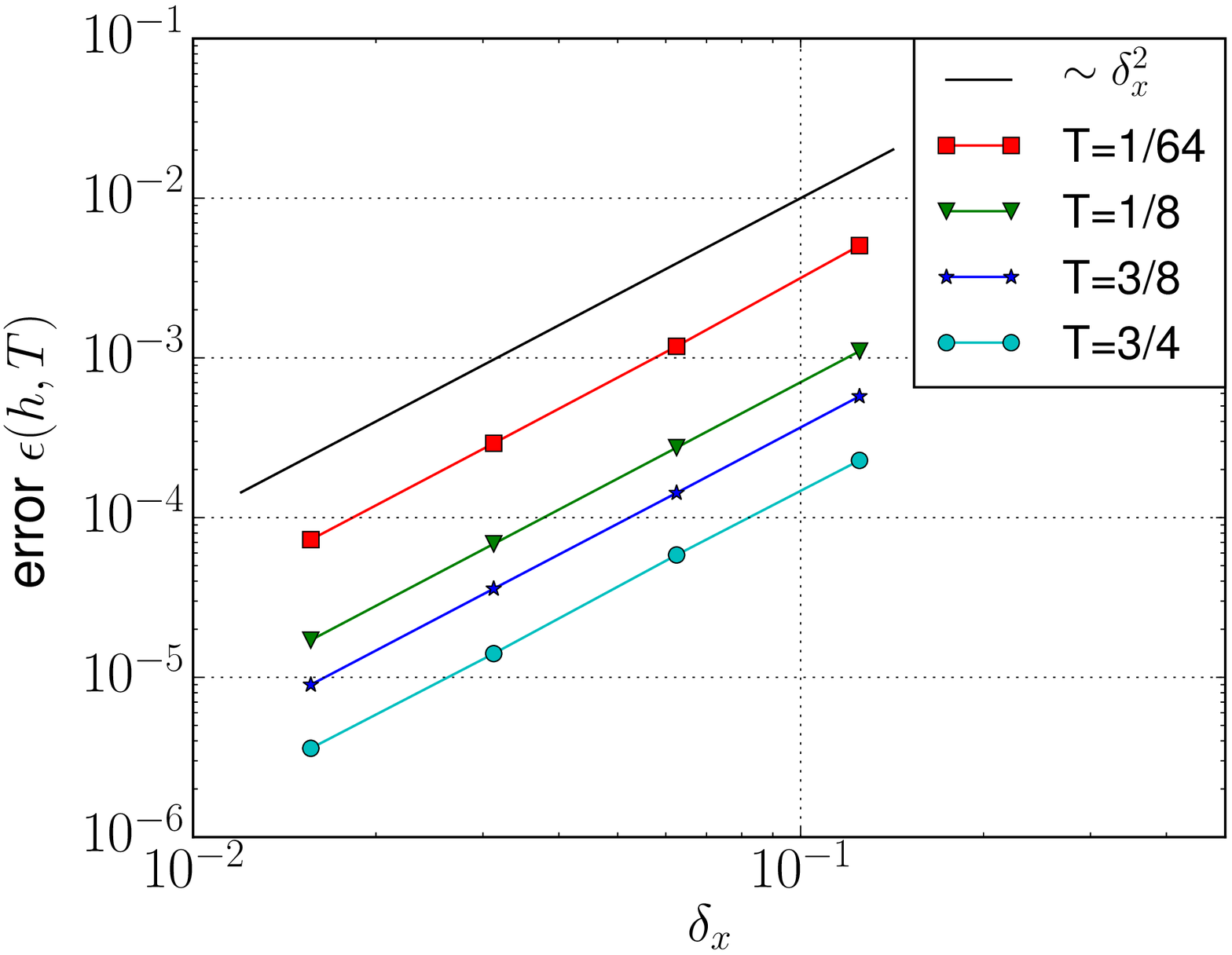}
  \includegraphics[width=0.49\linewidth]{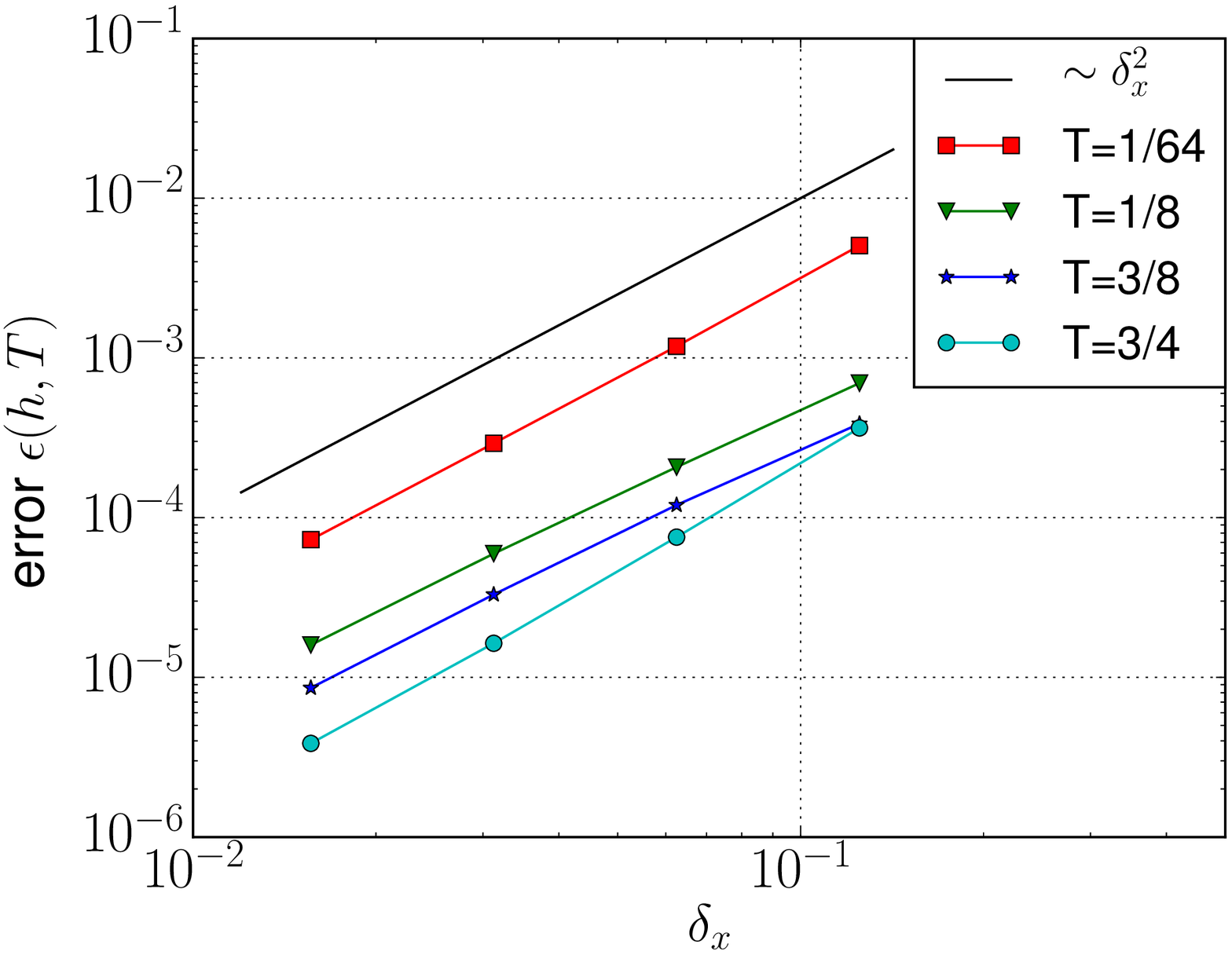}
  \caption{Grid convergence of the FSK-rule (left) and the new FSL-rule, based
    on linear interpolation (right) in plate driven flow at selected
    times. Since the channel width is an integral number, both approaches show a
    second order rate of convergence.}
  \label{fig:freeUniformConvergence}
\end{figure}

In our first validation case, we monitor the transient behavior of a planar flow with initial condition
$\vec{u}(\x,0) =0$. The domain is periodic in the $x$- and $y$-direction, with a
free surface boundary at $z=0$ and a solid wall at $z=h$, moving in the
$x$-direction with constant tangential velocity of $u_{wall} = 0.001$ lattice
units (Cf. Fig.~\ref{fig:sketch_channelFlow} with $\alpha=0$). 
This setup has been proposed by \citet{Yin2006} with the analytic Fourier series solution,
\begin{equation}
  \frac{u_{id}(z,t)}{u_{wall}} = 1 - 
  \sum_{k=0}^{\infty} { \frac{4 (-1)^k}{(2k+1)\pi} e^{-(2k+1)^2 \pi^2 \mu t / (4 \rho h^2)} \times \cos \left( \frac{(2k+1) \pi z}{2h} \right)  },
  \label{eq:freeUniform}
\end{equation}
for the validation of a free-slip boundary condition 
that leads in this case to the same solution as the free
surface condition.
For $t \gg 1$ the flow quickly develops into a uniform profile, since
the free surface does not impose any friction. A dimensionless time scale $T =
\mu t / (\rho h^2)$ is introduced to facilitate the evaluation of the flow at
the times $T=1/64$, $1/8$, $3/8$ and $3/4$. For the simulations, we use $\rho=1$
and $\mu=1/6$ in lattice units for channels of height $h=8$, $16$, $32$, $64$.
Qualitative results are shown in Fig.~\ref{fig:freeUniformProfile}. Note,
that very similar flow profiles are obtained for both the original free surface boundary
condition (FSK) and the newly proposed FSL-rule, since
here the channel height is restricted to have an integer value (in lattice cells).
For the quantitative evaluation, we define the error as
\begin{equation}
  \epsilon(h, T) = \frac{1}{u_{wall}} \sqrt{ \frac{1}{h} \sum_{z_i} (u_x(z_i,T) - u_{id}(z_i,T))^2 },
\end{equation}
where $z_i$ ranges over all the lattice node positions along the
$z$-axis. Fig.~\ref{fig:freeUniformConvergence} shows that both boundary
conditions yield correct transient behavior and the expected second order rate
of convergence is exhibited clearly.
The results have been obtained with a TRT - parametrization
of $\Lambda=1/4$.

\subsection{Linear Couette Flows}
The analysis predicts the exact recovery of linear flow profiles when the second
order boundary condition of Sec.~\ref{sec:shearbc} is used (FSL - rule with
prescribed boundary value $S_{\a \b}^b$). 
Here we evaluate  the case of a steady flow as follows. In a
cubic domain, we impose non-slip boundary conditions (bounce-back) at $z=0$,
fixing the position of the first lattice nodes to the plane $z=0.5$
(Cf. Fig.~\ref{fig:sketch_channelFlow} with $\alpha=0^\circ$). The shear rate
condition of Sec.~\ref{sec:shearbc} is imposed at $z=h$.  As a first
verification experiment, a tangential shear rate is imposed by setting $S^b_{x
  z}(z=h) = 0.001$. The steady Couette profile is recovered without numerical
error, independent of choice of equilibrium function and film thickness $h$, in
accordance with the analytical properties of the boundary condition.

\begin{figure}
  \centering
  \includegraphics[width=0.66\linewidth]{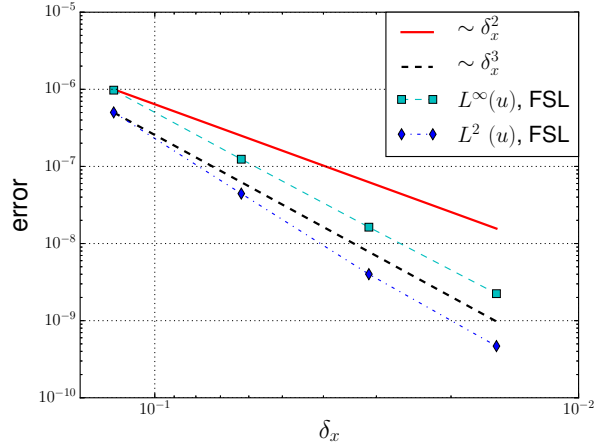}
  \caption{Second order convergence rate for linear shear flow when imposing
    constant shear rate on the top boundary in a 
    skewed channel with slope
    $\Delta_z / \Delta_x = 1/4$.}
  \label{fig:linearCouette}
\end{figure}

Our next test case is a rotated linear film flow where bottom and top boundary
planes are placed with a slope of $\Delta_z / \Delta_x = 1/4$ (i.e., $\alpha \approx 14^{\circ} $ in Fig.~\ref{fig:sketch_channelFlow}). In order to
realize the skew non-slip boundary, we use the CLI boundary condition of
\citet{Ginzburg2007}, which is a second order link-wise boundary condition
similar to the one proposed by \citet{Bouzidi2001}, based on linear
interpolation. This boundary condition can recover steady Couette flows in
arbitrary rotated channels exactly, provided that linear equilibria are
employed. Applying again a tangential shear rate $\partial_n u_{t} = 0.001$ the
exact profile is recovered if the equilibrium function is restricted
to the linear terms. If a non-linear equilibrium is used, a spurious
Knudsen-layer appears at the boundary nodes of the skew channel where the shear
rate is prescribed using Eq.~\eqref{eq:fsfit}.
From the analysis, we expect this error to be of second
order. A grid convergence study with fixed lattice viscosity $\nu$ and Reynolds
number $Re=0.064$ is conducted, varying channel widths $h_i = h_0, 2 h_0, 4 h_0,
8 h_0$ and imposed shear rate $\partial_n u_{t} = 0.001, 0.00025, 6.25e-05,
1.5625e-05$. Fig.~\ref{fig:linearCouette} shows that the grid convergence is
indeed of second order. Here and in the following sections, the relative errors
are computed using either the $L^2$-norm,
\begin{equation}
  L^2(\Phi) = \sqrt{ \frac{\sum_{\vec{x}}(\Phi(\vec{x}) - \Phi_{id}(\vec{x}))^2}{ \sum_{\vec{x}} \Phi_{id}(\vec{x})^2 }  },
  \label{eq:l2error}
\end{equation}
or the Tchebysheff norm,
\begin{equation}
  L^{\infty}(\Phi) = \frac{\max_{\vec{x}} \vert \Phi(\vec{x}) - \Phi_{id}(\vec{x}) \vert}{ \max_{\vec{x}} \vert \Phi_{id}(\vec{x}) \vert },
  \label{eq:maxerror}
\end{equation}
where $\Phi(\x)$ and $\Phi_{id}(\x)$ are the respective numerical and the ideal
value at the node position $\x$.

\subsection{Steady Parabolic Film Flow}
\label{sec:filmflow}
\begin{figure}
  \centering
  \includegraphics[width=0.66\linewidth]{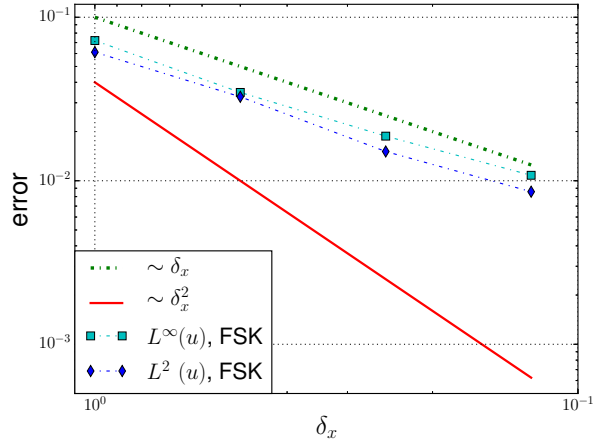}
  \caption{First order rate of convergence for a planar film flow of height
    $h=8.33$ with the FSK-rule. The same problem is solved exactly using
    the FSL-rule.}
  \label{fig:planarFilmFlowConvergenceFSK}
\end{figure}

Force-driven slow flow of finite thickness over a planar non-slip surface admits
an analytic solution that is used for validation as follows. Using a cubic
domain, we impose a non-slip boundary condition at the bottom $z=0$ plane of the
domain, realized using the bounce back rule. This means that the first lattice
nodes are located at a distance $0.5$ from the bottom plane. At $z=h$, a free
boundary is realized using the FSL boundary condition of Sec.~\ref{sec:shearbc}
with $S_{\a \b}^b=0$.
We use periodic boundary conditions in the $x$ and $y$ direction. The magic
parametrization $\Lambda = 3/16$ for parabolic straight channel flows is used
\citep{Ginzburg1994}, to eliminate the error of the bounce back
rule. 
It can be verified readily that the shear boundary condition yields the correct
steady state profile without numerical error, independent of the film thickness
$h$, and independent of choice of the equilibrium function. Applying additional
gravity directed towards the bottom plane yields an additional linear
hydrostatic pressure gradient that does not influence the solution, if
the ``incompressible equilibria'' \citep{HeLuo1997b} are used. 
Notice, that the FSK - rule of \citet{KoernerEtAl}
is exact in this test case only if $h$ is divisible by the grid spacing,
otherwise the expected accuracy is of first order
$\mathcal{O}(\delta_x)$. Fig.~\ref{fig:planarFilmFlowConvergenceFSK} shows that
the measured error convergence for the FSK - rule given by Eq.~\eqref{eq:koerner} is indeed reduced to first
order for $h=8.33$.

\begin{figure}
  \centering
  \includegraphics[width=0.66\linewidth]{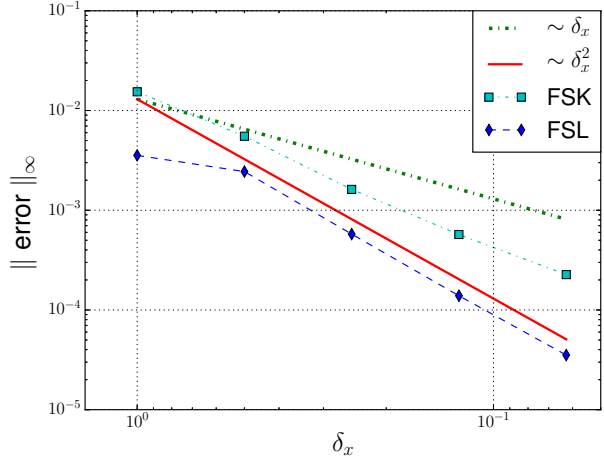}
  \caption{Comparison of the convergence behavior in a rotated planar film
    flow. The rate of convergence with the proposed FSL-rule is second order as
    predicted by the analysis. The behavior of the FSK-rule with $D=0$ is below
    second order.}
  \label{fig:planarFilmFlowSkewConvergence}
\end{figure}

We repeat the test case with the flow direction rotated about a slope of
$\Delta_z / \Delta_x = 1/7$ ($\alpha \approx 8.1^{\circ}$ in
Fig.~\ref{fig:sketch_channelFlow}) with respect to the lattice. The CLI boundary
condition is used for the skew non-slip wall to assure a second order rate of
convergence, and fix the parametrization using $\Lambda = 1/4$. Similar to 
Couette flow, now a certain error is inevitable. Using the interpolated FSL
boundary rule for the free boundary, an error convergence rate of order
$\mathcal{O}(\delta_x^2)$ is expected, independent of the flow direction,
opposed to a first order error for the original free surface condition from
Eq.~\eqref{eq:koerner} (FSK with $D=0$). The grid spacings are $\delta_x = 1,
0.5, 0.25, 0.125, 0.0625$, keeping the Reynolds number constant by adjusting the
accelerating force according to $g = g_0 \times \delta_x^{-3}$ at a constant
relaxation time $\tau = 2$. Fig.~\ref{fig:planarFilmFlowSkewConvergence} shows
the grid convergence of the two different boundary conditions. Indeed, the proposed 
FSL - condition shows a second order behavior, whereas for the original
FSK boundary condition the obtained rate of convergence is clearly below second
order. 


\subsection{Breaking Dam}
\begin{figure}
  \centering
  \includegraphics[width=0.66\linewidth]{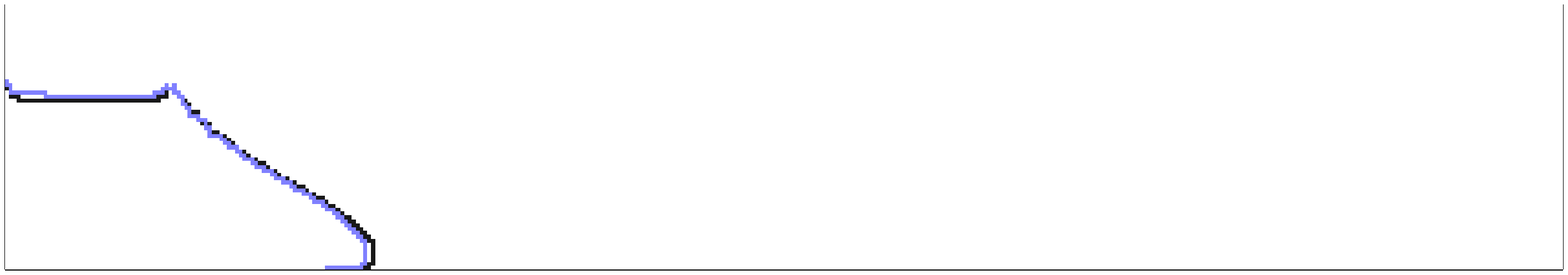}\\
  \includegraphics[width=0.66\linewidth]{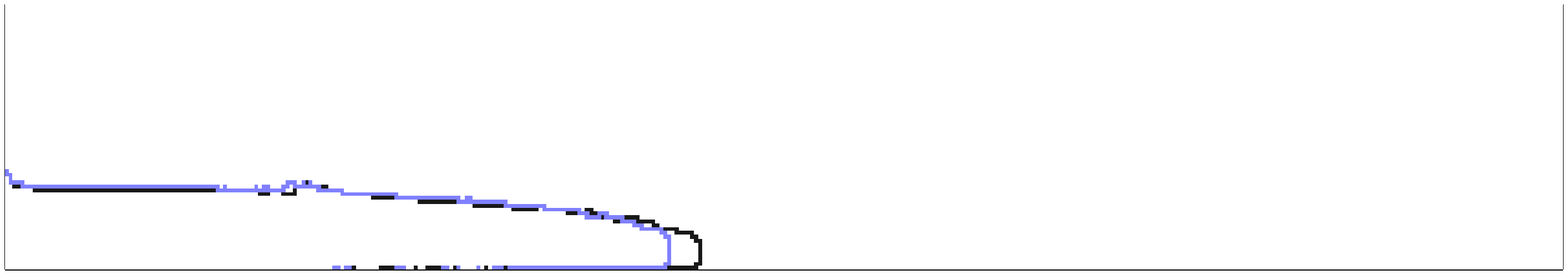}\\
  \includegraphics[width=0.66\linewidth]{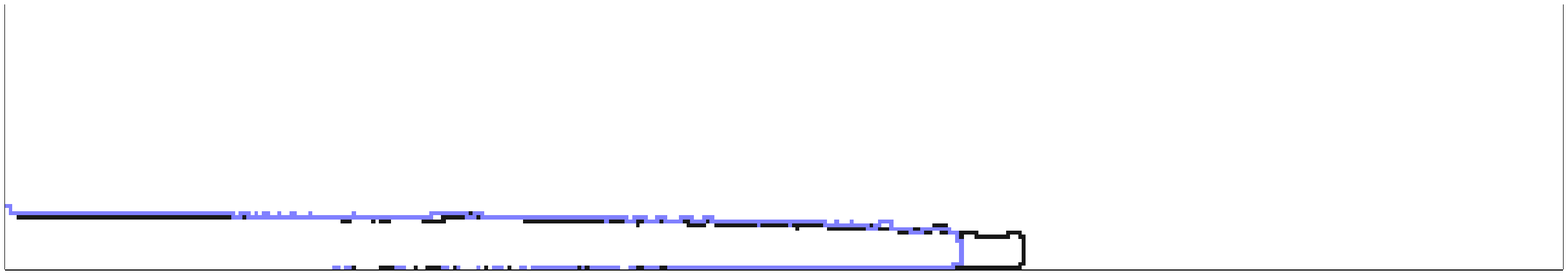}\\
  \includegraphics[width=0.66\linewidth]{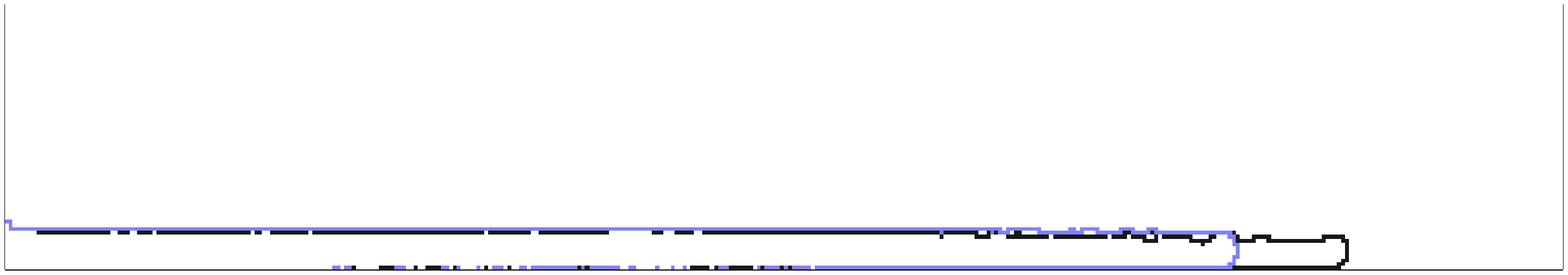}\\
  \caption{Simulation of a breaking dam after time step 2000, 4000, 6000, 8000.
    The faster surge front (black) is obtained from using the full boundary
    condition, while the slower front (gray blue) is obtained from the simplified
    boundary scheme with $D=0$ in Eq.~\eqref{eq:bterm}.}
  \label{fig:BreakingDam}
\end{figure}

Finally, we demonstrate the 
effect of the viscous stress term in the free
surface condition of Eq.s~\eqref{eq:proj1} and \eqref{eq:proj2} in the
instationary case of a collapsing rectangular column of liquid under gravity
(breaking dam). At the surge front, we have $\partial_n j_n > 0$.
However, the
simplified boundary rule with $D=0$ forces $\partial_n j_n = 0$ at the free
surface, hence we expect lower acceleration of the surge front for this case.
The simulated column has an initial size of $80$x$40$ lattice units. At a
lattice viscosity $\nu=1/3$ and a maximal flow velocity of $0.05$ this
corresponds to a Reynolds number of $Re = 12$.  Indeed,
Fig.~\ref{fig:BreakingDam} shows that the collapse of the column is
significantly slower if the terms are neglected in the boundary rule. For this
experiment, we 
use the first order FSK rule and only first order
(next-neighbor) extrapolation of $\rho$, $\vec{u}$ and $S_{\a \b}$ to compute
the boundary values. 
The interface tracking implementation (cf. Sec.~\ref{sec:fslbm}) is directly based on the original works \citep{KoernerEtAl,Thuerey2006}.

\section{\label{sec:conclusion}Conclusion}
Based on Chapman-Enskog analysis of the lattice Boltzmann equation we have
described the construction of boundary conditions for free surfaces with second
order spatial accuracy. In contrast to free surface models based on a 
discretization of the Navier-Stokes equations that need to impose boundary
values for the velocity, the free surface lattice Boltzmann approach imposes the
stress conditions directly on the distribution functions. Hence, the macroscopic
momentum 
appears in the closure relation only to 
match the non-linear terms.  The numerical experiments confirm the analytical
findings, i.e., that the proposed FSL boundary scheme is of second order spatial
accuracy, whereas the original FSK - model of \citep{KoernerEtAl} is only first
order accurate. However, in order to achieve full second order accuracy, the
interface position must be defined with the same order of accuracy to obtain the
correct $\delta$ - values of the link intersection with the boundary. This is
not possible with the interface tracking approach that is used in the original
FSLBM.  Hence, future implementations must make use of higher order interface
reconstruction methods, or other techniques such as level sets
\citep{Osher2001,Sethian2003} to represent the free surface. Inevitably, this
will introduce additional algorithmic complexity but eventually improve the
accuracy \citep{Nichols71}.

For full consistency with the defining equations of a free surface, the scheme
needs an approximation of the shear stress at the boundary, to impose the correct
boundary values on the LBM data. In a classical dam break problem at $Re=12$,
the significance of the viscous stresses at the boundary became visible.
At lower viscosities this term is probably less important. 
It should be noted, that for under-resolved free surface simulations, it is often more accurate to
employ a simplified boundary scheme, because physical viscosity and simulated viscosity
do not match. For instance, \citet{JanssenOffshore2010} have reported excellent
coincidence of high $Re$ - breaking dam simulations with experimental data using
the original FSK - rule neglecting the viscous terms. 
This effect in free surface simulations has already been described in \citep{Hirt68}.

The method developed in this article can serve as basis for 
a second order accurate LBM - based implementation of free surface flows when 
combined with second order accurate interface tracking.


\section*{Acknowledgements}
\noindent The authors would like to thank Irina Ginzburg for helpful discussion. 

The first author would like to thank the Bayerische Forschungsstiftung and
KONWIHR project waLBerla-EXA for financial support. The second author is
supported by the European Union Seventh Framework Program - Research for SME’s
with full title “High Productivity Electron Beam Melting Additive Man-
ufacturing Development for the Part Production Systems Market” and grant
agreement number 286695.

\bibliographystyle{plainnat}
\bibliography{lit}

\end{document}